\documentclass[notitlepage,11pt]{article}

\usepackage{amssymb,amsmath,comment,cases}

\usepackage{enumerate}
\usepackage{amsthm}

\usepackage{mathtools}
\usepackage[latin1]{inputenc}
\usepackage[english]{babel}

\theoremstyle{plain}
\newtheorem{theor}{Theorem}[]
\newtheorem{lemma}[theor]{Lemma}
\newtheorem{corolla}[theor]{Corollary}
\newtheorem{propos}[theor]{Proposition}
\newtheorem*{propos*}{Proposition}

\theoremstyle{definition}
\newtheorem{defin}{Definition}
\newtheorem{rem}[theor]{Remark}

\newcommand{\abs}[1]{\left\lvert {#1} \right\rvert}
\newcommand{\norma}[1]{\left\lVert {#1} \right\rVert}
\newcommand{\R}{\mathbb{R}}
\renewcommand{\phi}{\varphi}
\renewcommand{\epsilon}{\varepsilon}
\newcommand{\eps}{\varepsilon}
\renewcommand{\succ}[2]{\lbrace {#1}_{{#2}} \rbrace _{{#2} \in \NN}}

\newcommand{\RR}{\mathbb{R}}

\newcommand{\NN}{\mathbb{N}}

\def\f{\varphi}

\def\QED{\hfill {$\square$}\goodbreak \medskip}

\begin{document}

\title %{\vspace{-10mm}
{The non-anticoercive H\'enon-Lane-Emden system}

\author{Andrea Carioli\footnote{SISSA, via Bonomea, 265 -- 34136 Trieste, Italy. Email: {acarioli@sissa.it}.
{Partially supported by INDAM-GNAMPA.}}~ and
Roberta Musina\footnote{Dipartimento di Matematica ed Informatica, Universit\`a di Udine,
via delle Scienze, 206 -- 33100 Udine, Italy. Email: {roberta.musina@uniud.it}. 
{Partially supported by Miur-PRIN 201274FYK7\_004.}
}}

\date{}

\maketitle

\footnotesize

\noindent

\begin{abstract}{\small \noindent
We use variational methods to study the existence of a principal eigenvalue for the non-anticoercive
H\'enon-Lane-Emden system on a bounded domain. Then we provide a detailed insight into the problem in the
linear case.}
\end{abstract}

\normalsize
\vskip2.5cm

\section{Introduction}
The H\'enon-Lane-Emden system
\begin{equation}
\label{eq:HLE}
\begin{cases}
-\Delta u=|x|^a|v|^{p-2}v&\text{in $\Omega$}\\
-\Delta v=|x|^b|u|^{q-2}u&\text{in $\Omega$}\\
u=0=v&\text{on $\partial\Omega$}
\end{cases}
\end{equation}
includes the second and fourth order Lane-Emden equations
and the H\'enon equation in astrophysics. Here $\Omega$ is a 
domain in $\R^n$ containing the origin, and $a, b, p,q$ are given, with 
$a, b>-n$, $p,q>1$.

Most of the papers about problem (\ref{eq:HLE}) require $n\ge 3$ and
deal with the so-called \emph{anticoercive case}
$$
\frac{1}{p}+\frac{1}{q}<1.~\!
$$
Since the celebrated paper \cite{Lio1985} by P.L. Lions, where $a=b=0$ and $\Omega=\R^n$ are assumed,
large efforts have been made in investigating (\ref{eq:HLE})
and related problems. We will not try to provide a complete
list of references on this subject. We limit ourselves to cite
\cite{BVG, BMR1, BMR2, BusMan, CDM, DefigPerRos2008,
F, FG, Mit1, MS, PQS, SZ3, SZ2, SZ, Sou} and  references therein. 
Several challenging  problems are still open. Some of them
affect the  \emph{critical hyperbola} 
$$
\frac{a+n}{p}+\frac{b+n}{q}=n-2~\!,
$$
that was firstly introduced by Mitidieri \cite{Mit0} in 1990, in the autonomous case
$a=b=0$ (see also  \cite{CleDefigMit1992, PelVdvor1992, Mit1}). Roughly speaking,
to have existence of solutions one is lead to assume
\begin{equation}
\label{eq:existence}
\frac{a+n}{p}+\frac{b+n}{q}>n-2~\!.
\end{equation}

\medskip

The present paper deals with the case of a bounded and smooth domain $\Omega$ in $\R^n$, $n\ge 1$,
and with the {\em non-anticoercive case}
$$
q=p'=\frac{p}{p-1}~.
$$
Due to the homogeneities involved, we are lead to study
the  {\em eigenvalue problem}
\begin{equation}
\label{eq:problem}
\begin{cases}
-\Delta u=\lambda_1 |x|^a|v|^{p'-2}v&\text{in $\Omega$}\\
-\Delta v=\lambda_2 |x|^b|u|^{p-2}u&\text{in $\Omega$}\\
u=0=v&\text{on $\partial\Omega$~\!.}
\end{cases}
\end{equation}
We emphasize the fact that we include the lower dimensional cases
$n=1,2$, that actually present some peculiarities. 

As far as we know, only few references are available for (\ref{eq:problem}).
We mention the paper \cite{Mon2000}, where Montenegro uses
degree theory to face problem (\ref{eq:problem}) 
in a more general setting that includes non-self adjoint elliptic operators.

\medskip

We adopt a variational approach that allows us to weaken the 
integrability assumptions on the coefficients from Montenegro's $L^n(\Omega)$ to 
$L^1(\Omega)$. More precisely, we assume
\begin{subnumcases} {\label{eq:retta} }
~a,b>-n&{}\label{eq:retta_a}\\
~\frac{a}{p'}+\frac{b}{p}+2>0.&{}\label{eq:retta_b}
\end{subnumcases}
Notice that (\ref{eq:retta_b}) is automatically satisfied if $n=1, 2$ and (\ref{eq:retta_a}) holds.
Moreover,  (\ref{eq:retta_b}) coincides with  (\ref{eq:existence}), as $q=p'$ in our
setting. 

We
 look for \emph{finite energy solutions}, accordingly with the next definition.

\begin{defin}
\label{D:def} {\em
The pair $(u,v)$ is a finite-energy solution to (\ref{eq:problem}) if

\noindent
$\bullet$ $u,v\in W^{2,1}(\Omega)\cap W^{1,1}_0(\Omega)$;

\noindent
$\bullet$  $u\in L^p(\Omega, |x|^b dx)$, $v\in L^{p'}(\Omega, |x|^a dx)$, that is,
\begin{equation}\label{eq:finEn}
\int\limits_\Omega \abs{x}^b \abs{u}^p~\!dx < \infty~, \quad \int\limits_\Omega \abs{x}^a \abs{v}^{p'}~\!dx < \infty~\!;
\end{equation}
\noindent
$\bullet$ $u,v$ are weak solutions to the elliptic equations in 
(\ref{eq:problem}). That is,
$$
\int\limits_\Omega\!\nabla u\cdot\nabla \f~\!dx=\lambda_1\!\!\int\limits_\Omega|x|^a|v|^{p'-2}v\f~\!dx,
~~ \int\limits_\Omega\!\!\nabla v\cdot\nabla \f~\!dx=\lambda_2\!\!\int\limits_\Omega|x|^a|u|^{p-2}u\f~\!dx
$$
~~for any test function $\f\in C^\infty_c(\Omega)$.}
\end{defin}

Our approach is based on the equivalence\footnote{\footnotesize{already noticed for instance
 by Wang \cite{Wan1993} and Calanchi-Ruf \cite{CalRuf2010}  in the anticoercive case.}}
 between  \eqref{eq:problem} and 
the fourth order  eigenvalue problem
\begin{equation}
\label{eq:fourth}
\begin{cases}
\Delta\left( \abs{x}^{-a(p-1)} \abs{\Delta u}^{p-2}\Delta u \right) = \mu \abs{x}^b \abs{u}^{p-2}u&
\text{in $\Omega$}\\
u=\Delta u=0&\text{on $\partial\Omega$,}
\end{cases}
\end{equation}
where $\mu, \lambda_1$ and $\lambda_2$ satisfy
\begin{equation}\label{eq:Lambda}
\abs{\lambda_1}^{p-1}\lambda_1 \abs{\lambda_2}^{p'-1}\lambda_2 = \mu^{p'},
\end{equation}
compare with Lemma \ref{theor:equiv}.
We refer to Section \ref{S:approach} for the proof of the next result.

\begin{theor}\label{theor:main_intro}
Let $\Omega$ be a bounded and smooth domain in $\R^n$.
If (\ref{eq:retta}) hold, then problem (\ref{eq:problem})
has a positive principal eigenvalue $\mu$. That is, 
for any pair of real numbers $(\lambda_1, \lambda_2)$ satisfying
(\ref{eq:Lambda}),
problem (\ref{eq:problem}) has a finite-energy solution $(u,v)$,
such that $u,v>0$ in $\Omega$.
\end{theor}

\medskip
The last part of the paper is focused on the linear case  
$p=2$, so that (\ref{eq:retta_b}) becomes
\begin{equation}
\label{eq:retta2}
a+b+4>0~\!.
\end{equation}
In Section  \ref{S:linear} we prove that problem
\begin{equation}
\label{eq:p=2}
\begin{cases}
-\Delta u=\lambda_1 |x|^av&\text{in $\Omega$}\\
-\Delta v=\lambda_2 |x|^bu&\text{in $\Omega$}\\
u=0=v&\text{on $\partial\Omega$}
\end{cases}
\end{equation}
has a unique and simple {principal eigenvalue} $\mu_1>0$, and a discrete spectrum $\succ{\mu}{k}$. 
More precisely, the following results hold.

\begin{theor}
\label{T:p=2_1}
Let $\Omega$ be a bounded and smooth domain in $\R^n$. Let $a,b>-n$ and assume that
(\ref{eq:retta2}) holds. 
\begin{description}
\item$~i)$ There exists an increasing,
unbounded sequence of eigenvalues $\succ{\mu}{k}$ such that problem
(\ref{eq:p=2}) has a nontrivial and finite-energy solution $(u,v)$ if and only if $\lambda_1\lambda_2=\mu_k$
for some integer $k\ge 1$. 
\item $ii)$ The first eigenvalue $\mu_1$  is the unique principal eigenvalue. In addition,
$\mu_1$ is simple, that is, 
 if $(\tilde u,\tilde v)$ solves (\ref{eq:p=2}) and $\lambda_1\lambda_1=\mu_1$, then $\tilde u=\alpha u$ and $\tilde v=\beta v$
 for some $\alpha,\beta\in \R$.
 \end{description}
\end{theor}

\section{The variational approach}
\label{S:Preliminary}
In this section we introduce and study certain second order weighted Sobolev spaces
under Navier boundary conditions that are suitable for studying  (\ref{eq:fourth})
via variational methods.

\subsection{Functional setting and  embedding results}\label{SS:space}
To simplify notation we set $s=a(p-1)$. Thus, from now on we 
assume that $s, b$ are given exponents such that
$$
s>n-np~,\quad b>-n
$$
even if not explicitly stated. In addition, $\Omega \subset \RR^n$ will always denote
a bounded and smooth domain. We  denote by $c$ any universal positive constant.

We introduce the function space
$$
C^2_N(\overline{\Omega}) := \lbrace u \in C^2(\overline{\Omega}) ~\vert ~ u = 0 \text{ on } \partial\Omega
\rbrace.
$$
Let ${W}^{2, p}_N(\Omega, \abs{x}^{-s}dx)$ be the reflexive Banach space defined as the completion of the set
\[
{D}_0 := \lbrace u \in C^2_N(\overline{\Omega}) ~\vert~ \Delta u \equiv 0 \text{ on a neighborhood of
the origin}\rbrace~\!,
\]
with respect to the uniformly convex norm
\[
\norma{u}_{s} \equiv \norma{u}_{p, s} := \left( \int_\Omega \abs{x}^{-s} \abs{\Delta u}^p~\!dx \right)^\frac{1}{p}.
\]
We begin to study the spaces ${W}^{2, p}_N(\Omega,
\abs{x}^{-s}dx)$ by pointing out few embedding results. Firstly, notice that
the boundedness of the domain $\Omega$ implies
\begin{equation}\label{eq:scale}
{W}^{2, p}_N(\Omega, \abs{x}^{-{s}}~\!dx) \hookrightarrow {W}^{2, p}_N(\Omega,
\abs{x}^{-s_0}dx) \quad \text{ if } {s_0} \leq s.
\end{equation}
In order to simplify notation in the next lemma, we introduce the exponent 
$$
\hat p_s=p~~\text{if $s\ge 0$}~,\qquad
\hat p_s=\frac{np}{n-s}~~\text{if $s<0$}.
$$

\begin{lemma}\label{lem:car1}
Assume $s > n - np$. Then
$${W}^{2, p}_N(\Omega, \abs{x}^{-s}dx)\hookrightarrow W^{2, \tau}(\Omega)\cap W^{1, \tau}_0 (\Omega)
\quad\text{for any $\tau\in[1,\hat p_s)$.}
$$
\begin{proof}
Notice that $1<\hat p_s<p$.
For any $u\in D_0$  and $\tau\in [1,\hat p_s)$ we use elliptic regularity estimates, 
see for instance \cite[Lemma 9.17]{GilTru1983}, to get
$$
\|u\|^\tau_{W^{2,\tau}(\Omega)}\le c
 \int_\Omega \abs{\Delta u}^{\tau}~\!dx \leq c\left( \int_\Omega \abs{x}^{-s} \abs{\Delta u}^p ~\!dx
\right)^\frac{\tau}{p} \left( \int_\Omega \abs{x}^\frac{s\tau}{p-\tau}
~\!dx\right)^\frac{p-\tau}{p}.
$$
The last integral is finite as $s>n-np$, and the lemma readily is proved.
\end{proof}
\end{lemma}

The next lemma will be used in the next section to rigorously
prove the equivalence between the second order
system (\ref{eq:problem}) and the fourth order equation (\ref{eq:fourth}). 

\begin{lemma}\label{theor:caratt}
If $s>n-np$, then $u\in {W}^{2, p}_N(\Omega,\abs{x}^{-s}dx)$ if and only if 
\begin{equation}
\label{eq:cond}
u\in W^{2,1}\cap W^{1,1}_0(\Omega)~~\text{and}~~-\Delta u\in L^p(\Omega,|x|^{-s}dx).
\end{equation}
\begin{proof}
Clearly,  any  $u\in {W}^{2, p}_N(\Omega, \abs{x}^{-s}dx)$ 
satisfies (\ref{eq:cond}) by Lemma \ref{lem:car1}.

Conversely, fix $u$ satisfying (\ref{eq:cond}). Assume in addition that $-\Delta u = 0$ almost everywhere on a ball $B_r$ about $0$,
so that $-\Delta u\in L^p(\Omega)$. Hence, $u\in W^{2,p}(\Omega)\cap W^{1,p}_0(\Omega)$
by elliptic regularity theory. Extend $u$ to a function $u$ in $W^{2,p}(\R^n)$ with compact support
and take 
a sequence of mollifiers $\succ{\rho}{k}$. 
Since for $k$ large enough, $-\Delta(\rho_k*u)\equiv 0$ on $B_{r/2}$ and
$\rho_k*u\to u$ in $W^{2,p}(\Omega)$, then $-\Delta(\rho_k*u)\to -\Delta u$ in  $L^p(\Omega,|x|^{-s}dx)$. 
Let $u_k$ be the solution to
$$
\begin{cases}
-\Delta u_k=-\Delta(\rho_k*u)&\text{in $\Omega$}\\
u_k=0&\text{on $\partial\Omega$.}
\end{cases}
$$
It turns out that $u_k\in D_0\cap W^{2,p}_N(\Omega)$, as $u_k$ is smooth up to the boundary of $\Omega$
by regularity theory, and $-\Delta u_k\equiv 0$ in $B_{r/2}$. In addition, $u_k\to u$ in $W^{2,p}(\Omega)$
and $-\Delta u_k\to -\Delta u$  in  $L^p(\Omega,|x|^{-s}dx)$,
that is sufficient to conclude that $u\in  {W}^{2, p}_N(\Omega,\abs{x}^{-s}dx)$.

For a general $u$ satisfying (\ref{eq:cond}) let $u_k$ be the unique solution to 
\[\begin{cases}
-\Delta u_k = \chi_{\Omega_k}(-\Delta u)&\text{in $\Omega$}\\
u_k=0&\text{on $\partial\Omega$,}
\end{cases}
\]
where 
$\Omega_k := \Omega \setminus \overline{B}_{\eps_k}$ and $\eps_k\to 0$.
Then $u_k\in W^{2, p} \cap W^{1, p}_0(\Omega)$ and $u_k\in  {W}^{2, p}_N(\Omega,
\abs{x}^{-s}dx)$ by the first part of the proof. Clearly, the sequence $\succ{u}{k}$ is bounded in
${W}^{2, p}_N(\Omega, \abs{x}^{-s}dx)$, and we can assume that 
$u_k\to \bar u$ weakly in ${W}^{2, p}_N(\Omega, \abs{x}^{-s}dx)$. On the other hand,
$-\Delta u_k$ converges to $-\Delta u$ in $L^p(\Omega,  \abs{x}^{-s}dx)$ by Lebesgue's theorem.
Thus $\bar u=u$, that is, $u\in {W}^{2, p}_N(\Omega, \abs{x}^{-s}dx)$.
\end{proof}
\end{lemma}

The next corollary is an immediate consequence of  Lemma \ref{theor:caratt}.

\begin{corolla}\label{cor:caratt}
Assume $s>n-np$. For any $f \in L^p(\Omega,|x|^{-s}dx)$, the unique solution $u$ to
$$
\begin{cases}
 -\Delta u = f&\text{in $\Omega$}\\
u=0&\text{on $\partial \Omega$}
\end{cases}
$$
belongs to ${W}^{2, p}_N(\Omega, \abs{x}^{-s}dx)$.
\end{corolla}

Next we deal with embeddings in weighted $L^p$ spaces. 

\begin{lemma}\label{theor:immer}
If $s+b+2p\ge 0$, then 
\[
\Lambda(s, b) := \inf_{\substack{u \in {W}^{2, p}_N(\Omega, \abs{x}^{-s}dx)\\u \neq 0}}
 \frac{\displaystyle\int_\Omega \abs{x}^{-s} \abs{\Delta
u}^p~\!dx}{\displaystyle\int_\Omega{\abs{x}^b \abs{u}^p~\!dx}} > 0.
\]
\begin{proof}
First of all, notice that ${L}^{p}(\Omega, \abs{x}^{b_0}~\!dx) \hookrightarrow L^p(\Omega,
\abs{x}^{b}~\!dx)$ if ${b_0} \leq b$, that together with (\ref{eq:scale}) implies
\begin{equation}\label{eq:scaleLp}
\Lambda(s,b)\ge c \Lambda(s_0,b_0)\quad \text{if $s_0\le s$ and $b_0\le b$.}
\end{equation}
We start with the lowest dimensions $n=1,2$. Fix an exponent $s_0\le s$, such that
$n-np<s_0\le b(p-1)$. Then 
$\Lambda(s,b)\ge c \Lambda\big(s_0,\frac{s_0}{p-1}\big)>0$
by (\ref{eq:scaleLp}) and Lemma \ref{lem:bidIneq}
in the Appendix.

\medskip

Next, assume $n\ge 3$. In addition, assume firstly that $s< n-2p$. By a Rellich-type inequality in \cite{Mit2}, see also
{\cite[Lemma 2.14]{Mus2014}},  and using \cite[Lemma 2.9]{Mus2014}, one readily checks that
there exists a positive and explicitly known constant $c=c(n,p,s)$, such that 
\begin{equation}
\label{eq:Rellich}
c \int_\Omega \abs{x}^{-s-2p} \abs{u}^p~\!dx \leq \int_\Omega \abs{x}^{-s} \abs{\Delta u}^p~\!dx
\quad \text{for any $u \in \mathit{C}^2_N (\overline{\Omega})$,} 
\end{equation}
that is, $c=\Lambda(s,-s-2p)>0$. Thus $\Lambda(s,b)\ge c\Lambda(s,-s-2p)>0$ by (\ref{eq:scaleLp}). 
Finally, if $s\ge n-2p$, we fix a parameter $s_0$  such that 
$$
\max\{n-np,-2p-b\}<s_0<n-2p\le s,
$$
that is possible as $b>-n$ and $n\ge 3$. Then (\ref{eq:scaleLp}) and (\ref{eq:Rellich})  with $s$ replaced by $s_0$ give
$\Lambda(s,b)\ge c\Lambda(s_0,-s_0-2p)>0$, and the lemma is proved.
\end{proof}
\end{lemma}

\begin{rem}
{\em If $\Omega$ contains the origin and $s+b+2p<0$ then
$\Lambda(s,b)=0$. Indeed, fix a nontrivial $\psi \in C^\infty_c(B_1\setminus\{0\})$. 
For $k$ large enough the function $\psi_k(x)=\psi(kx)$
has compact support in $\Omega\setminus\{0\}$. Thus
\[
\Lambda(s, b) \le  \frac{\displaystyle\int_\Omega \abs{x}^{-s} \abs{\Delta \psi_k}^p ~\!dx}
{ \displaystyle\int_\Omega \abs{x}^b \abs{\psi_k}^p~\!dx}
= C k^{s+2p+b} = o(1) \quad \text{as $k \rightarrow \infty$}. 
\]}
\end{rem}

\begin{rem}\label{rem:notAppro}
{\em If $n - np<s< n-2p$, then $\mathit{C}^2_N (\overline{\Omega})\subset  {W}^{2, p}_N(\Omega, \abs{x}^{-s}dx)$
and the space
$$
\mathit{C}^2_N (\overline{\Omega} \setminus \lbrace 0 \rbrace) := \lbrace u \in
\mathit{C}^2_N(\overline{\Omega}) ~\vert~ u \equiv 0 \text{ on a neighborhood of the origin}\rbrace
$$
is dense in 
$ {W}^{2, p}_N(\Omega, \abs{x}^{-s}dx)$, see Lemma 2.14 in \cite{Mus2014}.
}
\end{rem}

\begin{rem}\label{rem:vuenne}
{\em By Lemma \ref{theor:caratt}, the set $D_0$ is dense in the standard
Sobolev space $W^{2,p}_N(\Omega)=W^{2,p}(\Omega)\cap W^{1,p}_0(\Omega)$.
The smaller set $C^2_N(\overline\Omega\setminus\{0\})$ is dense
in $W^{2,p}_N(\Omega)$ if and only if $n>2p$, compare with Remark
\ref{rem:notAppro}.}
\end{rem}

The next compactness result is a crucial point
for studying the
eigenvalue problem \eqref{eq:fourth}.

\begin{lemma}\label{theor:15}
If $s+b+2p>0$ then ${W}^{2, p}_N(\Omega, \abs{x}^{-s}dx)$ is compactly embedded into
$L^p(\Omega, \abs{x}^{b}dx)$.
\begin{proof}
It suffices to show that any sequence $\succ{u}{k}$ that converges weakly to the null function
in ${W}^{2, p}_N(\Omega, \abs{x}^{-s}dx)$ actually converges in $L^p(\Omega, \abs{x}^{b}dx)$.
Fix such a sequence, and take $\eps>0$ small. Since clearly $\succ{u}{k}$ is bounded in $W^{2,p}(\Omega
\setminus \overline B_\eps)$, then $|x|^b|u_k|^p\to 0$ in $L^1(\Omega\setminus \overline B_\eps)$
by Rellich theorem. Therefore, for any $b_0\in(-n,b)$ we have that
$$
\int_\Omega|x|^b|u_k|^p~\!dx=\int_{B_\eps}|x|^b|u_k|^p~\!dx+o(1)\le 
\eps^{b-b_0}\int_\Omega|x|^{b_0}|u_k|^p~\!dx+o(1)~\!.
$$
Now, if $b_0$ is close enough to $b$, then $s+b_0+2p>0$. Hence
$$
\int_\Omega|x|^b|u_k|^p~\!dx\le c\eps ^{b-b_0}+o(1)
$$
by Lemma \ref{theor:immer}.
The conclusion follows, as $\eps>0$ was arbitrarily chosen.
\end{proof}
\end{lemma}

\section{Two equivalent problems}\label{S:approach}
In this section we furnish a rigorous proof of the equivalence between the
eigenvalue problems (\ref{eq:problem}) and (\ref{eq:fourth}).
We start with a preliminary result.

\begin{lemma}\label{lemma:relation}
Assume that (\ref{eq:retta}) hold. For any $f \in L^p(\Omega, \abs{x}^{b}dx)$,
the  problem
\begin{subnumcases}{\label{eq:sistEffe}}
-\Delta u = \abs{x}^a \abs{v}^{p'-2}v&{}\label{eq:sistEffe_a}\\
-\Delta v = \abs{x}^b\abs{f}^{p-2}f&{}\label{eq:sistEffe_b}\\
\nonumber
u \in {W}^{2, p}_N(\Omega, \abs{x}^{-a(p-1)}~\!dx), ~v \in {W}^{2, p'}_N(\Omega,
\abs{x}^{-b(p'-1)}~\!dx)
\end{subnumcases}
admits a unique solution.
\begin{proof}
First of all, notice that $a(p-1) > n - np$, $b(p'-1) > n - np'$. Thus the results in 
Subsection \ref{SS:space} apply to the spaces ${W}^{2, p}_N(\Omega, \abs{x}^{-a(p-1)}~\!dx)$ and 
${W}^{2, p'}_N(\Omega,\abs{x}^{-b(p'-1)}~\!dx)$.

Since $\abs{x}^b \abs{f}^{p-2}f\in L^{p'}(\Omega, \abs{x}^{-b(p'-1)}~\!dx)$, then 
Corollary \ref{cor:caratt} guarantees (\ref{eq:sistEffe_b}) has a unique solution
$v\in {W}^{2, p'}_N(\Omega,\abs{x}^{-b(p'-1)}~\!dx)$. The embedding Lemma \ref{theor:immer} gives that
$\abs{x}^a\abs{v}^{p'-2}v\in L^p(\Omega, \abs{x}^{-a(p-1)}~\!dx)$.
Thus there exists a unique solution  $u\in {W}^{2, p}_N(\Omega, \abs{x}^{-a(p-1)}~\!dx)$ to (\ref{eq:sistEffe_a}), thanks
again to 
Corollary \ref{cor:caratt}.
\end{proof}
\end{lemma}

We are ready to prove the claimed equivalence result.
\begin{lemma}\label{theor:equiv}
Assume that  (\ref{eq:retta}) hold. Let $\mu, \lambda_1,\lambda_2\in \R$ satisfying 
(\ref{eq:Lambda}). Then the following sentences are equivalent.
\begin{description}
\item$~~i)$ $u \in {W}^{2, p}_N(\Omega, \abs{x}^{-a(p-1)}~\!dx)$ is a weak solution to (\ref{eq:fourth}).
\item$~ii)$ The pair $u,v:= -|x|^{-a(p-1)}|\Delta u|^{p-2}\Delta u$ solves
\begin{subnumcases}{\label{eq:sisCrit}}
\label{eq:sisCrit_a}
-\Delta u = \abs{x}^a \abs{v}^{p'-2}v\\
\label{eq:sisCrit_b}
-\Delta v = \mu \abs{x}^b\abs{u}^{p-2}u,\\
\nonumber
u\in {W}^{2, p}_N(\Omega, \abs{x}^{-a(p-1)}~\!dx)~,~v\in {W}^{2,p'}_N(\Omega, \abs{x}^{-b(p'-1)}~\!dx).
\end{subnumcases}
\item$iii)$ The pair $u,v:= -|x|^{-a(p-1)}|\Delta u|^{p-2}\Delta u$ is a finite-energy solution to 
(\ref{eq:problem}), in the sense of Definition \ref{D:def}.
\end{description}
\begin{proof}
If $u \in {W}^{2, p}_N(\Omega, \abs{x}^{-a(p-1)}~\!dx)$, then $u\in L^p(\Omega,|x|^b dx)$
by Lemma \ref{theor:immer}. Thus we can apply Lemma \ref{lemma:relation} to find a unique pair 
$u_0, v_0$ such that
\[
\begin{cases}
-\Delta u_0 = \abs{x}^a \abs{v_0}^{p'-2}v_0\\
-\Delta v_0 = \mu \abs{x}^b\abs{u}^{p-2}u,\\
u_0 \in {W}^{2, p}_N(\Omega, \abs{x}^{-a(p-1)}~\!dx), v_0 \in {W}^{2,
p'}(\Omega, \abs{x}^{-b(p'-1)}~\!dx).
\end{cases}
\]
Notice that
$v_0=|x|^{-a(p-1)}|\Delta u_0|^{p-2}(-\Delta u_0)$ almost everywhere in $\Omega$.
Therefore, if $u$ solves (\ref{eq:fourth}), then for any $\f\in D_0$ it holds that
\begin{eqnarray*}
\int_\Omega \abs{x}^{-a(p-1)} \abs{\Delta u}^{p-2} \Delta u \Delta \f ~\!dx &=&
\mu\int_\Omega|x|^b\abs{u}^{p-2} u \f ~\!dx\\
&=&\int_\Omega(-\Delta v_0)\f~\!dx=\int_\Omega v_0(-\Delta\f)~\!dx\\
&=& \int_\Omega \abs{x}^{-a(p-1)} \abs{\Delta u_0}^{p-2} \Delta u_0 \Delta \f ~\!dx,
\end{eqnarray*}
that readily gives that $u=u_0$, since $u,u_0\in  {W}^{2, p}_N(\Omega, \abs{x}^{-a(p-1)}~\!dx)$
and $D_0$ is dense in ${W}^{2, p}_N(\Omega, \abs{x}^{-a(p-1)}~\!dx)$. Hence, also $v=v_0$, the pair $u,v$ solves (\ref{eq:sisCrit}),
and the first implication is proved.

The equivalence between $ii)$ and $iii)$ is immediate, thanks to 
Lemma \ref{lem:car1} and Corollary  \ref{cor:caratt}. It remains to
show that $ii)$ implies $i)$.
If $(u, v)$ solves \eqref{eq:sisCrit}, then for every $\f \in D_0$ it holds that
\begin{eqnarray*}
\mu \int_\Omega \abs{x}^b \abs{u}^{p-2}u \f ~\!dx &=& 
\int_\Omega v (-\Delta \f )~\!dx \\
&=&
\int_\Omega \abs{x}^{-a(p-1)} \abs{\Delta u}^{p-2} \Delta u \Delta \f ~\!dx,
\end{eqnarray*}
that is, $u$ solves (\ref{eq:fourth}).
\end{proof}
\end{lemma}
Lemma \ref{theor:equiv}  shows that finite energy solutions to 
(\ref{eq:problem}) are the stationary points of the functional
$$
u\mapsto \int_{\Omega}\abs{x}^{-a(p-1)}|\Delta u|^p~\!dx
$$
on the constraint 
$$
M = \left\{ u\in {W}^{2, p}_N(\Omega, \abs{x}^{-a(p-1)}~\!dx) ~\bigg\vert~ \int_\Omega|x|^b|u|^p~\!dx=1~\right\}~\!.
$$
If (\ref{eq:retta}) hold, then $M$ is compact in the weak 
${W}^{2, p}_N(\Omega, \abs{x}^{-a(p-1)}~\!dx)$ topology, by Lemma \ref{theor:15}. Thus
the infimum 
\begin{equation}
\label{eq:mu}
\mu:= \Lambda(a(p-1), b) = \inf_{\substack{u\in {W}^{2, p}_N(\Omega,
\abs{x}^{-a(p-1)}~\!dx) \\ u \neq 0}} \frac{\displaystyle\int_\Omega \abs{x}^{-a(p-1)} \abs{\Delta
u}^p~\!dx}{\displaystyle\int_\Omega{\abs{x}^b \abs{u}^p~\!dx}}
\end{equation}
is positive and attained. The next lemma deals with minimizers for $\mu$.

\begin{lemma}
\label{L:u_pos}
Assume that (\ref{eq:retta}) hold. 
If $u\in {W}^{2, p}_N(\Omega,\abs{x}^{-a(p-1)}~\!dx)$ achieves  
$\mu$, then, up to a change of sign, $u$ is positive
and superharmonic.
\begin{proof}
Let $v=|x|^{-a(p-1)}|\Delta u|^{p-2}(-\Delta u)$, so that the pair $(u,v)$ solves
 \eqref{eq:sisCrit}. Use Corollary \ref{cor:caratt} to introduce $u_0$ via
\[
\begin{cases}
- \Delta u_0 = \abs{x}^a \abs{v}^{p' - 1}\\
u_0 \in {W}^{2, p}_N(\Omega, \abs{x}^{-a(p-1)}~\!dx).
\end{cases}
\]
In particular, $u_0$ is superharmonic and
positive on $\Omega$. Next, put 
\[
g = \abs{x}^a \abs{v}^{p'-2}v.
\]
Thus $u$ and $u_0$ solve, for some $\tau\in[1,p)$,
\[
\begin{cases}
-\Delta u = g\\
u \in W^{2, \tau} \cap W^{1, \tau}_0 (\Omega),
\end{cases}
\qquad
\begin{cases}
-\Delta u_0 = \abs{g}\\
u_0 \in W^{2, \tau} \cap W^{1, \tau}_0 (\Omega).
\end{cases}
\]
Since $-\Delta(u_0 \pm u ) \geq 0$ and $u_0\pm u=0$ on the boundary of $\Omega$,
then $u_0 \pm u \geq 0$, that is, $u_0 \geq \abs{u}$. On the other hand, $\abs{\Delta
u_0} = \abs{g} = \abs{\Delta u}$. Therefore 
$$
\mu\le \frac{\displaystyle\int_\Omega \abs{x}^{-a(p-1)} \abs{\Delta u_0}^p~\!dx}
{\displaystyle\int_\Omega{\abs{x}^b \abs{u_0}^p~\!dx}} \le
\frac{\displaystyle\int_\Omega \abs{x}^{-a(p-1)} \abs{\Delta u}^p~\!dx}{\displaystyle\int_\Omega{\abs{x}^b \abs{u}^p~\!dx}} =
\mu,
$$
that is, $u_0$ attains $\mu$ and $u_0 = \abs{u}$. Since $u_0$ is positive in $\Omega$, then $u$ and
$-\Delta u$ have constant sign, as desired.
\end{proof}
\end{lemma}

\medskip
\noindent
{\bf Proof of Theorem \ref{theor:main_intro}.}
Immediate, thanks to Lemmata \ref{theor:15} and \ref{L:u_pos} and  the
equivalence given by Lemma \ref{theor:equiv}.
\QED

We conclude the section pointing out a symmetry result about the infimum
in (\ref{eq:mu}). It is convenient to use the notation $\mu(a,b,p)$ to 
emphasize the dependence of $\mu$ on the parameters.
\begin{propos}\label{theor:symmetry}
If  (\ref{eq:retta}) hold, then
$\mu(b, a, p')^p = \mu(a, b, p)^{p'}$.

\begin{proof}
Let $u$ be an extremal for $\mu(a,b,p)$. By Lemma \ref{theor:equiv}, the pair 
$u, v=-|x|^{-a(p-1)}|\Delta u|^{p-2}\Delta u$
solves \eqref{eq:sisCrit} with $\mu=\mu(a,b,p)$ and hence
\begin{multline*}
\int_\Omega \abs{x}^{-b(p'-1)} \abs{\Delta v}^{p'}~\!dx = 
\mu(a, b, p)^{p'} \int_\Omega \abs{x}^{b} \abs{u}^{p}~\!dx\\
\quad = \mu(a, b, p)^{p'-1} \int_\Omega \abs{x}^{-a(p-1)} \abs{\Delta u}^p~\!dx = 
\mu(a, b, p)^{p'-1} \int_\Omega \abs{x}^a \abs{v}^{p'}~\!dx.
\end{multline*}
Thus $\mu(b, a, p') \leq \mu(a, b, p)^{p'-1}$, or
equivalently $\mu(b, a, p')^p \leq \mu(a, b, p)^{p'}$. Exchanging the roles of $u$ and $v$, we get the opposite
inequality.
\end{proof}
\end{propos}

\section{Proof of Theorem \ref{T:p=2_1}}
\label{S:linear}
Since $p=2$, then equation (\ref{eq:fourth}) reduces to 
\begin{equation}\label{eq:equazEffe}
\Delta \left( \abs{x}^{-a} \Delta u\right) = \mu \abs{x}^b u.
\end{equation}

\bigskip
\noindent
{\em Proof of $i)$.}
We
denote by $X_a$
the Hilbert space ${W}^{2, 2}_N(\Omega, \abs{x}^{-a}dx)$, endowed with norm
$\norma{\cdot}_a$  and scalar product $(~\cdot~ | ~\cdot~)_a$.

We formally introduce  the``solution operator'' to \eqref{eq:equazEffe} under Navier boundary conditions.
More precisely, we define the linear operator $${T} \colon L^2(\Omega, \abs{x}^{b}dx) \rightarrow X_a~,
\quad
({T} f | w)_a = \int_\Omega \abs{x}^b f w~\!dx\quad \text{for every }w \in X_a.
$$
Then $T$ is  continuous, positive and self-adjoint. Let $j\colon X_a \rightarrow L^2(\Omega, \abs{x}^{b}dx)$ 
be the embedding in Lemma \ref{theor:immer}.
Then the operator 
$$T \equiv j \circ {T}~,~T \colon L^2(\Omega, \abs{x}^{b}dx) \rightarrow L^2(\Omega,
\abs{x}^{b}dx)$$ is compact. Thus the point spectrum $\sigma_p(T)$ of $T$ 
is a sequence $\succ{\nu}{k}$ of positive numbers converging
to 0, and
$$\displaystyle{
 \frac{1}{\nu_k} = \min \left\lbrace \frac{\displaystyle\int_\Omega \abs{x}^{-a} \abs{\Delta u}^2~\!dx}
{\displaystyle\int_{\Omega} \abs{x}^b u^2 ~\!dx} \colon u \in \Lambda_i^\perp, 1\leq i \leq k-1 \right\rbrace},
$$
where $\Lambda_i := \Lambda(\nu_i)$ is the eigenspace relative to the eigenvalue $\nu_i$.

By the results in the previous sections, $u$ is an eigenfunction for $T$ if and only
if the couple $(u,-\abs{x}^{-a}\Delta u)$
is the only solution to \eqref{eq:sistEffe}. This concludes the proof.

\bigskip
\noindent
{\em Proof of $ii)$.}
We will use the theory of abstract positive operators on Banach lattices, for which we refer to the
monograph~\cite{Schae1974}. Recall that $L^2(\Omega, \abs{x}^{b}dx)$
has a natural Banach lattice structure induced by the cone $P_+$ of nonnegative functions. 
We will show that $T$ is positive and irreducible. Then, the conclusion will follows thanks to 
an adaptation of Theorem V.5.2 in \cite{Schae1974}, that guarantees that the following facts hold:
\begin{enumerate}[(a)]
 \item The spectral radius $r(T) \in \RR_+$ is an eigenvalue.
\item The eigenspace $\Lambda(r(T))$ has dimension one, and is spanned by a (unique, normalized)
quasi-interior point of $P_+$.
\item $r(T)$ is the unique eigenvalue of $T$ with a positive eigenvector.
\end{enumerate}

To check that $T$ is irreducible we first recall that that the only closed ideals in $L^2(\Omega,
\abs{x}^{b}~\!dx)$ are the ones of the form
\[
 I_A = \left\lbrace f \in L^2(\Omega,
\abs{x}^{b}dx) \Big\vert f = 0 \text{ on } A \right\rbrace,
\]
where $A$ is a measurable set, see for instance \cite[p. 157]{Schae1974}. 
Therefore, we have to show that if $A$ satisfies
\[
 0 < \int_\Omega \abs{x}^{b} \chi_A ~\!dx < \int_\Omega \abs{x}^{b}~\!dx,
\]
then $I_A$ is not fixed by $T$. 

Let $f\in I_A$ be a nonnegative fixed function.
Then the problem
\[
\begin{cases}
-\Delta v = \abs{x}^b f\\
v \in X_b
\end{cases}
\]
admits a solution by Corollary \ref{cor:caratt}, and $v \in W^{2, \tau} \cap W^{1, \tau}_0(\Omega)$ for some $\tau > 1$. The minimum principles imply that $v$ is strictly positive in $\Omega$. For the same reason, the problem
\[
\begin{cases}
-\Delta u = \abs{x}^a v\\
u \in X_a,
\end{cases}
\]
defines a function $u$ that is strictly positive in $\Omega$. Hence $u \equiv Tf \notin I_A$, and this
proves the irreducibility property.
The same argument proves also the positivity property.
\QED

\section*{Appendix. An inequality in lower dimensions}

We sketch here the proof of some second order integral estimates in low dimensions by using, in
essence, the Rellich-type identity in \cite{Mit2}. Details and further applications of the underlying ideas will be given in
\cite{Car}.

\begin{lemma}\label{lem:bidIneq}
Assume $n=1$ or $2$ and let $\Omega$ be an open interval or 
a bounded domain in $\R^2$ of class  $C^2$.
If $s> n-np$, then there exists a constant $c>0$ such that
\[
c \int_\Omega \abs{x}^{\frac{s}{p-1}} \abs{u}^p ~\!dx \le \int_\Omega \abs{x}^{-s} \abs{\Delta u}^p
~\!dx
\]
for any $u \in C^2_N(\overline\Omega)$ such that $\Delta u=0$ in a neighborhood of $0$.
\begin{proof}
We can assume that $\Omega$ is contained in the unit ball about the origin. Put
$$
a=\frac{s}{p-1},
$$
and notice that $a>-n$. 
We argue in a heuristic way.
A more rigorous proof requires a suitable approximation of the weight $|x|^{a+2}$ by smooth
functions, see also \cite{Car}. 

Fix $u\in C^2_N(\overline\Omega)$ such that $\Delta u=0$ in a neighborhood of $0$. For $p\ge 2$ one clearly has
$$
(p-1)\int_\Omega  \abs{\nabla u}^2 \abs{u}^{p-2}~\!dx=\int_\Omega(-\Delta u)|u|^{p-2}u~\!dx.
$$
For general $p>1$, one can check that
$\abs{\nabla u}^2 \abs{u}^{p-2}\in L^1(\Omega)$ and
\begin{equation}\label{eq:div2}
(p-1)\int_\Omega  \abs{\nabla u}^2 \abs{u}^{p-2}~\!dx
\le
\int_\Omega  |\Delta u|\abs{u}^{p-1}~\!dx ~\!.
\end{equation}
Next, we are allowed to use integration by parts again and   H\"older inequality
 to estimate
 \begin{multline*}
(a+2)(a+n)\int_\Omega|x|^a|u|^p~\!dx=-
\int_\Omega (\Delta|x|^{a+2}) \abs{u}^p~\!dx \\=
\!p\int_\Omega (\nabla |x|^{a+2}  \cdot \nabla u) \abs{u}^{p-2} u~\!dx
\leq\!p(a+2)\int_\Omega |x|^{a+1}|\nabla u||u|^{p-1}~\!dx\\
\leq\!p(a+2)
\left( \int_\Omega\abs{\nabla u}^2 \abs{u}^{p-2}dx\right)^{\!\frac{1}{2}}
\left(\int_\Omega |x|^{2a+2}\abs{u}^p~\!dx\right)^{\!\!\frac{1}{2}}.
\end{multline*}
Therefore, from $a+2\ge a+n>0$ we infer
$$
\left(\frac{a+n}{p}\right)^2\int_\Omega |x|^a  \abs{u}^p~\!dx\le \int_\Omega  \abs{\nabla u}^2
\abs{u}^{p-2}~\!dx\le  \frac{1}{p-1}\int_\Omega  \abs{u}^{p-1} |\Delta u|~\!dx
$$
by \eqref{eq:div2}. 
It remains to use H\"older inequality to estimate
$$
c \int_\Omega |x|^a \abs{u}^p~\!dx\le \left( \int_\Omega |x|^a \abs{u}^p~\!dx \right)^\frac{1}{p'} \left( \int_\Omega
|x|^{-a(p-1)} \abs{\Delta u}^p \!dx \right)^\frac{1}{p}
$$
where $c=c(a,n,p)>0$. Thus 
$$
c \int_\Omega |x|^a \abs{u}^p~\!dx\le \int_\Omega
|x|^{-a(p-1)} \abs{\Delta u}^p \!dx,
$$
as desired. 
\end{proof}
\end{lemma}

\end{document}